\documentclass[11pt,fleqn]{article}
\usepackage{amsfonts}
\usepackage{mathrsfs}
\usepackage{amsthm}
 \usepackage{graphicx} 
\usepackage{amsmath} 
\usepackage{amssymb} 
\usepackage{enumerate}
\usepackage{epsfig}
\usepackage{multirow}
\usepackage{color}
\usepackage{cite}

\usepackage{epsfig}

\newtheorem{conj}{\hskip\parindent Conjecture}
\newtheorem{qe}{\hskip\parindent Question}
\newtheorem{thm}{\hskip\parindent Theorem}
\newtheorem{lem}{\hskip\parindent Lemma}
\newtheorem{prop}{\hskip\parindent Proposition}
\newtheorem{cor}{\hskip\parindent Corollary}

\numberwithin{equation}{section}

\def\Re {\mathop{\rm Re}\nolimits}
\def\Im {\mathop{\rm Im}\nolimits}

\begin{document}

\title{Proof of a conjecture of Granath on  optimal bounds  of the  Landau constants }
\author{Chun-Ru Zhao,  Wen-Gao Long and  Yu-Qiu Zhao\footnote{Corresponding author ({\it{E-mail address:}} {stszyq@mail.sysu.edu.cn}).  Investigation  supported in part by the National
Natural Science Foundation of China under grant numbers  10871212 and  11571375.}  }
  \date{ {\it{Department of Mathematics, Sun Yat-sen University, GuangZhou
510275, China}}
}

\maketitle

\begin{abstract}
We study the asymptotic expansion for the Landau constants $G_n$,
\begin{equation*}
\pi G_{n}\sim \ln(16N)+\gamma+\sum^{\infty}_{k=1}\frac{\alpha_k}{N^k} ~~\mbox{as} ~  n\rightarrow\infty,
\end{equation*}
where  $N=n+1$, and $\gamma$ is Euler's constant. We   show that the  signs of  the coefficients $\alpha_{k}$  demonstrate a periodic behavior such that
$(-1)^{\frac {l(l+1)} 2} \alpha_{l+1}< 0$ for  all $l$. We  further prove a conjecture of Granath which states that
$(-1)^{\frac {l(l+1)} 2} \varepsilon_l(N)<0$ for $l=0,1,2,\cdots$ and $n=0,1,2,\cdots$,  $\varepsilon_l(N)$ being the error due to truncation at the $l$-th order term.
 Consequently, we  also obtain the sharp bounds up to arbitrary orders of the form
 \begin{equation*}
 \ln(16N)+\gamma+\sum_{k=1}^{p}\frac{\alpha_{k}}{N^{k}}<\pi G_{n}<\ln(16N)+\gamma+\sum_{k=1}^{q}\frac{\alpha_{k}}{N^{k}}
 \end{equation*}
for all $n=0,1,2\cdots$,  all  $p=4s+1,\; 4s+2$ and $q=4m,\; 4m+3$, with  $s=0,1,2,\cdots$ and  $m=0, 1, 2,\cdots$.
\end{abstract}

\vspace{5mm}

 \vskip.3cm

\section{Introduction and Statement of Results}
\setcounter{equation} {0}

In 1913, Landau \cite{Landau} proved that
if $f(z)$ is analytic in the unit disc, and $\left|f(z)\right| <1$ for $|z|<1$, with the Maclaurin expansion
\begin{equation*}
f(z)=a_0+a_1 z+a_2 z^2+\cdots+a_n z^n+\cdots,~~ \left|z\right|<1,
\end{equation*}
then    there exist constants $G_n$ such that
\begin{equation*}
 \left| a_0+a_1+\cdots+a_n \right|\leq G_n,~~n=0,1,2,\cdots,
\end{equation*} and the   bound is optimal for each $n$,
  where $G_0=1$, and
\begin{equation}\label{Landau-constants}
G_{n}=1+\left (\frac{1}{2}\right )^{2}+\left (\frac{1\cdot3}{2\cdot4}\right )^{2}+\cdots
+\left (\frac{1\cdot3\cdot\cdots(2n-1)}{2\cdot4\cdot\cdots(2n)}\right )^{2}~~\mbox{for}~~n=1,2,\cdots.
\end{equation}

The constants $G_n$ are termed the Landau constants.
The large-$n$ behavior is known from the very beginning. Landau \cite{Landau}  derived  that
\begin{equation*}
G_n\sim \frac{1}{\pi}\ln n ~~\mbox{as}~n\rightarrow\infty;
\end{equation*}see also Watson \cite{Watson}.
It is worth mentioning that
 there exist  generating functions for these constants; cf. \cite{Cvijovic-Srivastava},  possible $q$-versions of the constants; cf. \cite{ismail-li-rahman2015}, and
an observation   made by Ramannujan (cf. \cite{Cvijovic-Srivastava}) that  relates   the Landau constants to the generalized hypergeometric functions.
Useful integral representations for $G_n$ have been   obtained  from such relations; cf., e.g.,   Watson \cite{Watson}; see also Cvijovi\'{c} and
Srivastava \cite{Cvijovic-Srivastava}.

The approximation of $G_n$ has gone in two related directions.  One is to obtain large-$n$ asymptotic approximations for the constants,  in  a   time period spanning from the early    twentieth century \cite{Landau, Watson}
to very recently  \cite{Cvijovic-Srivastava, li-liu-xu-zhao1}.
The other direction is   to find sharper bounds of $G_n$ for all nonnegative  integers $n$.
Authors working on the sharper bounds includes Brutman \cite{Brutmam} and Falaleev \cite{Falaleev} (in terms of elementary  functions),
Alzer \cite{Alzer} and  Cvijovi\'c and   Klinowski \cite{Cvijovic-Klinowski} (using the digamma function),
Zhao \cite{zhao},   Mortici \cite{Mortici} and Granath \cite{Granath} (involving higher order terms), and Chen and Choi \cite{chen-choi2014} and Chen \cite{chen2014} (digamma function and higher order terms).   The list is by no means complete.    The reader is referred to  \cite{Cvijovic-Srivastava, li-liu-xu-zhao1,li-liu-xu-zhao2} for a historic account.

\subsection{Optimal bounds up to all orders}
Attempts have been made to seek  bounds in a sense optimal, and up to arbitrary accuracy.

In 2012,   Nemes  \cite{nemes}  derived  full asymptotic expansions. For  $0 <h < 3/2$, he shows that the Landau constants $G_n$ have the   asymptotic expansion
\begin{equation}\label{nemes-expansion}
G_n\sim \frac 1 \pi \ln (n+h) +\frac 1 \pi (\gamma+4\ln 2 ) - \sum_{k\geq 1}\frac {g_k(h)}{(n+h)^k}~~\mbox{as}~~n\rightarrow +\infty,\end{equation}  where $\gamma=0.577215\cdots$ is   Euler's constant.
Earlier in 2011, the special cases $h=\frac 12$ and $h=1$ were established by  Nemes and Nemes    \cite{nemes-nemes}  using a formula in \cite{Cvijovic-Klinowski}.
  They also conjecture in \cite{nemes-nemes} a symmetry  property of  the computable constant coefficients
such that $g_k(h)=(-1)^k g_k(3/2-h)$ for every $k\geq 1$.
The conjecture has been  proved by G. Nemes himself   in  \cite{nemes}. A natural consequence is that for $h=3/4$,  all odd terms in the expansion vanish.
In this  important special case, Nemes \cite{nemes} has further proved that

\noindent
\begin{prop}
(Nemes) The following asymptotic approximation holds:
\begin{equation}\label{beta-alternative}
\pi G_n\sim   \ln (n+3/4) + \gamma+4\ln 2  + \sum_{s=1}^\infty \frac { \beta_{2s}}{ (n+3/4)^{2s}},~~n\rightarrow \infty, \end{equation}
where the coefficients $(-1)^{s+1}\beta_{2s}$ are positive rational numbers.\end{prop}

The derivation of Nemes \cite{nemes} is based on an integral representation of $G_n$ involving a  Gauss hypergeometric function in the integrand. An entirely different difference equation approach is applied in Li {\it et al.} \cite{li-liu-xu-zhao1} to obtain full asymptotic expansions with coefficients iteratively given. What is more, in a follow-up paper \cite{li-liu-xu-zhao2}, it is shown that the error due to truncation of \eqref{beta-alternative}
is bounded in absolutely value by, and of the same sign as, the first neglected term for all $n=0,1,2,\cdots$. An immediate   corollary is

\noindent
\begin{prop}
(Li, Liu, Xu and Zhao)
For $N=n+3/4$,
it holds
\begin{equation}\label{beta-bounds}
                   \ln N+\gamma+4\ln 2+\sum_{s=1}^{2m}\frac{\beta_{2s}}{N^{2s}}<     \pi G_n <  \ln N+\gamma+4\ln 2+\sum_{s=1}^{2k-1}\frac{\beta_{2s}}{N^{2s}}
 \end{equation}
 for all $n=0,1,2,\cdots$, $m=0,1,2,\cdots$, and $k=1,2,\cdots$.
 \end{prop}

In a sense, the formulas  \eqref{beta-alternative} and \eqref{beta-bounds} in the above propositions seem to have ended  a journey since one has thus obtained optimal bounds up to arbitrary orders. Yet there is an interesting observation worth mentioning,
as presented in the 2012 paper  \cite{Granath} of Granath; see also \cite{li-liu-xu-zhao2}.

Granath derives an asymptotic expansion
\begin{equation}\label{Granath-expansion}
\pi G_{n}\sim  \ln(16N)+\gamma +\sum_{k=1} ^\infty  \frac{\alpha_k}{N^k}, ~~n\rightarrow\infty,   \end{equation}  where $\alpha_k$ are  effectively computable  constants   but not explicitly given, except for the first few.
Here and hereafter we use the notation $N=n+1$.

Denoting the truncation of \eqref{Granath-expansion}
\begin{equation}\label{Granath-truncated}
A_l(N)=\ln(16N)+\gamma +\sum_{k=1} ^l  \frac{\alpha_k}{N^k},  \end{equation}
then one of the main results in Zhao \cite{zhao} reads $A_2(N) < \pi G_n< A_3(N)$ for $n=0,1,2,\cdots$.  Mortici \cite{Mortici} have actually proved that
 $
A_5(N)<\pi G_{n}<A_4(N)$ for   all non-negative $n$.

In \cite{Granath}, Granath proves that
$A_5(N) <\pi G_{n} < A_7(N)$ and states that   $A_9(N) <\pi G_{n} < A_{11}(N)$,   for all non-negative $n$.
Based on these formulas  and numerical evidences, Granath proposes a conjecture.

\noindent \begin{conj}(Granath)  It holds
 \begin{equation}\label{Granath-conjecture}
(-1)^{\frac {l(l+1)} 2} \left (  \pi G_{n}-A_l(N)\right ) <0   \end{equation}for all  $n=0, 1,2,\cdots$ and $l=0,1,2,\cdots$.
\end{conj}

\subsection{Statement of results}
We will show that the conjecture is true. To do so, we will make use of the second order difference equation for $G_n$ employed in \cite{li-liu-xu-zhao1}, and some estimating techniques used in \cite{li-liu-xu-zhao2}.

First we denote the error term
\begin{equation}\label{remainder}
\varepsilon_l(N)=\pi G_n -A_l(N)=\pi G_n -\left\{ \ln(16N)+\gamma +\sum_{k=1} ^l  \frac{\alpha_k}{N^k} \right\};
\end{equation}cf. \eqref{Granath-truncated}, where $N=n+1$. It is readily seen that $\varepsilon_l(N)\sim  {\alpha_{l+1}}/{N^{l+1}}$ as $N\to\infty$.
Hence we may  start by showing that \eqref{Granath-conjecture} holds for large $n$. To this aim, we have
\begin{thm}\label{coefficient}
The coefficients of the asymptotic expansion \eqref{Granath-expansion}  satisfy
\begin{equation}\label{the-coefficient-sgn}
 (-1)^{\frac {l(l+1)} 2} \alpha_{l+1}< 0,~~l=0,1,2,\cdots.
\end{equation}
\end{thm}

Next, we will prove the conjecture for all non-negative $n$.
\begin{thm} \label{thm epsilon-N}
For $N=n+1$,  it holds
\begin{equation}\label{the-remainder-sgn}
(-1)^{\frac {l(l+1)} 2} \varepsilon_l(N)<0
\end{equation}
for $l=0, 1,2,\cdots$ and  $n=0, 1,2,\cdots$.
\end{thm}

As a straightforward application of Theorem \ref{thm epsilon-N}, we obtain the following  sharp bounds up to arbitrary orders.
\begin{cor}\label{thm-inequalities}
For $N=n+1$, it holds $A_p(N)  <\pi G_{n}<A_q(N)$, that is,
\begin{equation}\label{sharp-bounds}
\ln(16N)+\gamma+\sum_{k=1}^{p}\frac{\alpha_{k}}{N^{k}}   <\pi G_{n}<   \ln(16N)+\gamma+\sum_{k=1}^{q}\frac{\alpha_{k}}{N^{k}}
\end{equation}
for all $n=0,1,2,\cdots$ and for all $p=4s+1,\;4s+2$ and  $q=4m,\;4m+3$, with $s=0,1,2,\cdots$ and $m=0,1,2,\cdots$.
\end{cor}
In view of Theorem \ref{coefficient}, we see that the bounds in  \eqref{sharp-bounds} are optimal as $n\to\infty$.

Theorem \ref{thm epsilon-N} can actually   be understood as an estimate of the error term, such that
 the error due to truncation is bounded in absolute value by, and of the same sign as, the first one or two neglected terms.
Indeed, since
    \begin{equation*}
  \varepsilon_l(N)=\frac {\alpha_{l+1}}{N^{l+1}}+ \varepsilon_{l+1}(N)~~\mbox{and}~~\varepsilon_l(N)=\frac {\alpha_{l+1}} {N^{l+1}}+ \frac { \alpha_{l+2}}{N^{l+2}}+ \varepsilon_{l+2}(N),\end{equation*}
taking into account the signs in Theorems \ref{coefficient} and \ref{thm epsilon-N}, we have
\begin{equation*}
  0< \varepsilon_{4k+1}(N)< \frac { \alpha_{4k+2}}{N^{4k+2}}+ \frac { \alpha_{4k+3}}{N^{4k+3}}~~\mbox{and}~~0< \varepsilon_{4k+2}(N)<   \frac { \alpha_{4k+3}}{N^{4k+3}}
\end{equation*}for all non-negative integers  $n$ and $k$, and
\begin{equation*}
 \frac { \alpha_{4k+1}}{N^{4k+1}} < \varepsilon_{4k}(N)<0~~\mbox{and}~~            \frac { \alpha_{4k+4}}{N^{4k+4}} +   \frac { \alpha_{4k+5}}{N^{4k+5}}  < \varepsilon_{4k+3}(N)< 0
\end{equation*}for all non-negative integers  $n$ and $k$.

As a by-product of the proof of Theorem \ref{thm epsilon-N}, we have  approximations of the asymptotic coefficients, follows respectively from  \eqref{alpha-2k-tuta-approx} and \eqref{alpha-tidle-2k+1}:
\noindent\begin{cor}\label{cor-coefficient-approx}
Assume that $\alpha_k$ are the coefficients in the asymptotic expansion \eqref{Granath-expansion}. Then we have
\begin{equation}\label{alpha-even-approx}
 \alpha_{2k}=(-1)^{k+1}  \frac { 2(2k-2)!} {(2\pi)^{2k}} \left (1+O\left (\frac 1 k\right )\right )
\end{equation}and
\begin{equation}\label{alpha-odd-approx}
\alpha_{2k+1}=(-1)^{k+1} \frac {8 (2k)! \ln (2k+1)}{(2\pi)^{2k+2}} \left ( 1+ O \left (\frac 1 {\ln k}\right )\right )
\end{equation} as $k\to \infty$.
\hfil\qed
\end{cor}

\section{The asymptotic coefficients and the proof of Theorem \ref{coefficient}}

From the representation \eqref{Landau-constants} one obtains the recurrence relation
\begin{equation*}
  G_{n+1}-G_n=  \left (\frac {2n+1}{2n+2}\right )^2 \left (G_n-G_{n-1}\right ).
\end{equation*}
Set $N=n+1$, we may rewrite it as a standard   second-order difference equation
\begin{equation} \label{difference-equation}
w(N+1)-\left (2-\frac  1 N+\frac 1{4N^2}\right ) w(N)+\left (1-\frac 1 {2N}\right )^2 w(N-1)=0,
\end{equation}
where $w(N)=\pi G_n$. An interesting fact is that the formal solution to \eqref{difference-equation} is an asymptotic solution; cf. Li and Wong \cite{Wong-Li1992a}; see also \cite{li-liu-xu-zhao1}.
Hence the asymptotic series  \eqref{Granath-expansion} furnishes a formal solution of \eqref{difference-equation}.
Therefore, one way to determine the coefficients $\alpha_k$ is to substitute \eqref{Granath-expansion} into \eqref{difference-equation} and equalizing  the coefficients of the same powers of $x=1/N$. We include some details as follows.
\begin{equation*}
 \ln (1+x)+\sum^\infty_{k=1} \frac {\alpha_k x^k}{(1+x)^k}
-\left (2-x+\frac {x^2} 4\right ) \sum^\infty_{k=1} \alpha_k x^k
+\left (1-\frac x 2\right )^2 \left [ \ln (1-x)
+\sum^\infty_{k=1}\frac{\alpha_k x^k}{(1-x)^k}\right ]=0.
\end{equation*}
Using  the Maclaurin series expansions,   we have
\begin{equation*}
 - \sum_{s=3}^\infty  d_{0,s} x^s +\sum^\infty_{k=1} \alpha_k x^k \sum^\infty_{j=2} d_{k, j+k}x^j=
  \sum^\infty_{s=3} \left ( \sum^{s-2}_{k=1} d_{k, s} \alpha_k - d_{0, s}\right ) x^s =0.
\end{equation*}
Accordingly,  coefficients $\alpha_k$ are determined  by
\begin{equation}\label{coefficient-alpha}
 d_{s-2, s}\alpha_{s-2}+d_{s-3,s}\alpha_{s-3}+\cdots+d_{1, s}\alpha_1-d_{0,s}=0,   ~~s=3,4, \cdots,
\end{equation}
where the coefficients $d_{s-2, s}=(s-2)^2$ for $s=3,4,\cdots$,
\begin{align}
                   & d_{0,s}=\frac {(-1)^{s} +1} s -\frac 1{s-1}+ \frac 1 {4(s-2)}~~\mbox{for}~s=3,4,\cdots,~~\mbox{and}\label{d-0-s} \\
                  & d_{k,s}=\frac {\left ( (-1)^{s-k} +1\right ) (k)_{s-k}} {(s-k)!}  -\frac {(k)_{s-k-1}}{(s-k-1)!}+ \frac {(k)_{s-k-2}} {4(s-k-2)!}    \label{d-k-s}
               \end{align}
for $k=1,2,\cdots, s-3$ and $s=k+3, k+4,\cdots$.

Appealing to
  \eqref{coefficient-alpha}-\eqref{d-k-s}, the first few coefficients $\alpha_k$ can be  evaluated as
\begin{equation*}
 \begin{array}{llll}
\alpha_1= -\frac 1 4 ,             &\alpha_2=\frac{5}{192} ,                       &\alpha_3=\frac{3}{128}   ,                             &\alpha_4= -\frac{341}{122880},         \\[.1cm]
\alpha_5= -\frac{75}{8192},        &\alpha_6=\frac{7615}{8257536}   ,              & \alpha_7=\frac{2079}{262144}  ,                       &\alpha_8= -\frac {679901}{1006632960} ,\\[.1cm]
\alpha_9=-\frac{409875}{33554432}, &\alpha_{10}=\frac{16210165}{17716740096}  ,    & \alpha_{11}=\frac{31709469}{1073741824}      ,        & \alpha_{12}=-\frac{568756771963}{281406257233920} .
   \end{array}
\end{equation*}
One readily sees   a periodic phenomenon of the signs of the coefficients, which agrees with Theorem \ref{coefficient}. To give a full proof  of the theorem, we may connect the coefficients with those in \eqref{beta-alternative}, and eventually with a certain hypergeometric function.
Indeed, re-expanding the formula \eqref{beta-alternative} in descending powers of  $N=n+1$ yields the  expansion    \eqref{Granath-expansion}. Hence we have
\begin{equation}\label{alpha-beta}
\alpha_k=4^{-k} \left [ -\frac 1 k+\sum^k_{j=1} \frac {(k-1)! 4^j \beta_j} {(j-1)! (k-j)!}\right ],~~k=1,2,\cdots;
\end{equation}
cf. \cite[(4.4)]{li-liu-xu-zhao2}, where $\beta_j$ vanish for odd integers $j$.
We also note that the coefficients $\beta_{2k}$ possess a   generating function, that is,
\begin{equation}\label{u-function}
u(x)=  \sum^\infty_{k=0} \rho_k  x^{2k}   =\frac x {2\sin\frac x 2} F\left ( \frac 1 2, \frac 1 2; 1; \sin^2\frac x 4\right ):=  \frac x {2\sin\frac x 2} F\left (  \sin^2\frac x 4\right )  ,
\end{equation}
where $\rho_0=1$ and $\rho_s=\frac {(-1)^{s+1}\beta_{2s}}{(2s-1)!}$, $s=1,2,\cdots$ are the  positive constants  defined in \cite[Sec. 3.1]{li-liu-xu-zhao2}.
 It is shown  in \cite{li-liu-xu-zhao2} that the generating function $u$ solves a second-order differential equation, and consequently  the hypergeometric function $F\left ( \frac 1 2, \frac 1 2; 1; t\right )$
 is brought in. It is worth noting that the function  also furnishes   a generating relation for the Landau constants, namely
$\frac {F(x)}{1-x}=\sum^\infty_{n=0} G_n x^n$ for small $x$; see \cite{nemes}.
 Here and hereafter we denote for short the hypergeometric function  as $F(t)=F\left ( \frac 1 2, \frac 1 2; 1; t\right )$.

\vskip .4cm

\noindent
{\bf{Proof of Theorem \ref{coefficient}}}. From \eqref{alpha-beta} we have
\begin{equation}\label{alpha-2k}
  \frac {\alpha_{2k}}{(2k-1)!} =\sum^k_{s=0} \frac {(-1)^{s+1} \rho_s}{(2k-2s)!} \left (\frac 1 4\right )^{2k-2s}~~\mbox{for}~k=1,2,\cdots .
\end{equation}Here use has been made of  the fact that $\beta_{2s-1}=0$ for $s=1,2,\cdots$. From \eqref{alpha-2k} we further have
\begin{equation}\label{generating-function-2k1}
  1+\sum^\infty_{k=1} \frac {(-1)^{k+1} \alpha_{2k}}{(2k-1)!}x^{2k} =\left\{\sum^\infty_{s=0}   \rho_s x^{2s}  \right\} \left\{
  \sum^\infty_{s=0}  \frac 1 {(2s)!}  \left (-\frac {x^2} {16} \right )^{s}\right\} =u(x)\cos\frac x 4 .
\end{equation}
Combining  \eqref{u-function}  with \eqref{generating-function-2k1}, and applying   a quadratic transformation formula,  we have
\begin{equation*}
 1+\sum^\infty_{k=1} \frac {(-1)^{k+1} \alpha_{2k}}{(2k-1)!}x^{2k} =\frac {\frac x 4}{\sin\frac x 4} F\left ( \sin^2\frac x 4\right )=  \frac {\frac x 4}{\sin\frac x 4} \frac 1 {\cos^2\frac x 8} F\left ( \tan^4\frac x 8\right );
\end{equation*}see \cite[(15.3.17)]{as}.
Each factor on the right-hand side   possesses  a Maclaurin expansion with positive coefficients; see Nemes \cite[pp.\;842-843]{nemes}. Hence we conclude that
\begin{equation}\label{alpha-even-sgn}
(-1)^{k+1} \alpha_{2k}>0~~\mbox{for}~k=1,2,\cdots.
\end{equation}

Similarly, we may write
\begin{equation}\label{tlide a odd}
  \sum^\infty_{k=0} \frac {(-1)^{k+1} \alpha_{2k+1}}{(2k)!}x^{2k}=\frac 1 4 \left\{\sum^\infty_{s=0}   \rho_s x^{2s}  \right\} \left\{
  \sum^\infty_{s=0}  \frac 1 {(2s+1)!}  \left (-\frac {x^2} {16} \right )^{s}\right\} =u(x)\frac {\sin\frac x 4} x .
\end{equation}Taking \eqref{u-function} into account, we can write the right-hand side term as
\begin{equation*}
  \frac 1 {4\cos\frac x 4}  F\left (   \sin^2\frac x 4\right )= \frac 1 {4\cos\frac x 4}   \frac 1 {\cos^2\frac x 8} F\left (  \tan^4\frac x 8\right ),
\end{equation*}which again has a Maclaurin expansion with all positive coefficients. Here we have used the formula
\begin{equation*}
  \frac 1 {\cos t}=\sum^\infty_{k=0}\frac {(-1)^k E_{2k} }{(2k)!} t^{2k},
\end{equation*}where $E_{2k}$ are the Euler numbers  such that $(-1)^k E_{2k}>0$ for $k=0,1,2,\cdots$; see \cite[(24.2.6)-(24.2.7)]{nist}.
Accordingly we have
\begin{equation}\label{alpha-odd-sgn}
 (-1)^{k+1} \alpha_{2k+1}>0~~\mbox{for}~k=0,1,2,\cdots.
\end{equation}
A combination of \eqref{alpha-even-sgn} and \eqref{alpha-odd-sgn} then gives
\eqref{the-coefficient-sgn}. \hfil\qed\vskip.5cm

\section{Proof of Theorem \ref{thm epsilon-N}}

To give a rigorous proof of Theorem \ref{thm epsilon-N}, we introduce
  \begin{equation}\label{R-l}
R_l(N)=\varepsilon_l(N+1)-\left (2-\frac 1 N +\frac 1 {4N^2}\right )     \varepsilon_l(N)+ \left (1-\frac 1 {2N}\right )^2  \varepsilon_l(N-1)
\end{equation}for $l=0,1,2\cdots$ and  $N=n+1=1,2,3,\cdots$, where $\varepsilon_l$ is the remainder term given in \eqref{remainder}.
 Similar to the derivation of \eqref{coefficient-alpha}, substituting \eqref{remainder} into  \eqref{R-l}, and again denoting $x=1/N$, we see that $R_l(N)$ is an analytic function of $x$ at the origin, with the Maclaurin expansion
\begin{equation}\label{R-l-expansion}
 R_l(N)=   \sum_{k=3}^\infty  d_{0,k} x^k -\sum^l_{k=1} \alpha_k x^k \sum^\infty_{j=2} d_{k, j+k}x^j  =\sum^\infty_{s=l+3} r_{l,s} x^s,
\end{equation}
where, for $s=l+3, l+4,\cdots$, and $l=0,1,2,\cdots$, the coefficients in \eqref{R-l-expansion} are
\begin{equation}\label{R-l-coefficient}
r_{l,s}=-\left(d_{l, s}\alpha_{l}+d_{l-1, s}\alpha_{l-1}+\cdots+d_{1, s}\alpha_1-d_{0, s}\right).
\end{equation}

To justify Theorem \ref{thm epsilon-N}, we state a lemma as follows,   leaving  the proof of it to later sections.
\begin{lem}{\label{lemma-R(N)}}For $N=n+1$, it holds
\begin{equation}\label{R-even-sgn}
 \tilde R_{2l}(N):=(-1)^{l+1} R_{2l}(N)>0,~~n=1,2,3,\cdots,~~l=0,1,2,\cdots .
\end{equation}
\end{lem}
Now we prove the theorem, assuming that  Lemma \ref{lemma-R(N)} holds true.\vskip .5cm

\noindent
{\bf{Proof of Theorem \ref{thm epsilon-N}}}. For fixed $l$, $l=0,1,2\cdots$,
first  we show that
\begin{equation}\label{remainder-even-sgn}
  \tilde\varepsilon_{2l}(N):=(-1)^{l+1} \varepsilon_{2l}(N)>0
\end{equation}
for  all  $n=1,2,\cdots$, where $N=n+1$. To this aim,
we note that
   \begin{equation}\label{tilde-varepsilon-2l}
  \tilde\varepsilon_{2l}(N) = (-1)^{l+1} \varepsilon_{2l}(N)= \frac {(-1)^{l+1}\alpha_{2l+1}}{N^{2l+1}}\left\{ 1+O\left ( \frac 1  N\right )\right \}
  =   \frac {\left |\alpha_{2l+1}\right |}{N^{2l+1}}\left\{ 1+O\left (  \frac 1  N\right )\right \}
  >0
  \end{equation}
 for $N$ large enough; cf. \eqref{Granath-expansion},  \eqref{remainder} and \eqref{alpha-odd-sgn}. Now assume that  $\tilde \varepsilon_{2l}(N)>0$ is not true for some $N$. Then there exists a finite positive $M$ such that
\begin{equation*}
  M=\max\{N=n+1: n\in \mathbb{N}~ \mbox{and}~ \tilde\varepsilon_{2l}(N)\leq 0\}.
\end{equation*}Thus for the positive integer $M$, we have   $\tilde\varepsilon_{2l}(M)\leq 0$, while $\tilde\varepsilon_{2l}(M+1),~ \tilde\varepsilon_{2l}(M+2),~ \cdots > 0$.

Denoting $b(N)= \left (1-\frac 1 {2N}\right )^2$ for simplicity, from \eqref{R-l} we have
\begin{equation*}
\tilde \varepsilon_{2l}(M+2)=(1+b(M+1))  \tilde \varepsilon_{2l}(M+1)+b(M+1) ( -\tilde \varepsilon_{2l}(M)) +\tilde R_{2l}(M+1).
\end{equation*}The later terms on the right-hand side are nonnegative  (where $M+1\geq 2$), hence we have
\begin{equation*}
\tilde \varepsilon_{2l}(M+2)\geq (1+b(M+1))  \tilde \varepsilon_{2l}(M+1) >\tilde \varepsilon_{2l}(M+1).
\end{equation*}
Moreover, from \eqref{R-even-sgn} we further have
\begin{equation*}
\tilde \varepsilon_{2l}(M+3)\geq (1+b(M+2))  \tilde \varepsilon_{2l}(M+2)+b(M+2) ( -\tilde \varepsilon_{2l}(M+1))
>  \tilde \varepsilon_{2l}(M+2) .\end{equation*}
Repeating  the process gives
\begin{equation*}
\tilde \varepsilon_{2l}(M+k+1)
>  \tilde \varepsilon_{2l}(M+k),~~k=1,2,\cdots  .\end{equation*}
By induction we conclude
\begin{equation}\label{inequality}
 \tilde \varepsilon_{2l}(M+1) < \tilde \varepsilon_{2l}(M+k)
\end{equation}for $k\geq 2$. Recalling that $\tilde \varepsilon_{2l}(N)=O\left ( N^{-2l-1}\right )$ for $N\to\infty$; cf. \eqref{tilde-varepsilon-2l}, letting $k\to\infty$ in \eqref{inequality} gives $\tilde \varepsilon_{2l}(M+1)\leq 0$. This contradicts the fact   that $\tilde\varepsilon_{2l}(M+1)>0$. Hence \eqref{remainder-even-sgn} holds.

 Now from    \eqref{remainder}, \eqref{alpha-even-sgn} and \eqref{remainder-even-sgn},   we  have
 \begin{equation}\label{remainder-odd-sgn}
 (-1)^l \varepsilon_{2l-1}(N) = \frac {  (-1)^{l} \alpha_{2l}} {N^{2l}} +(-1)^l \varepsilon_{2l}(N)<0~~\mbox{for}~~ l=1,2,3,\cdots,~~\mbox{and}~~ n=1,2,\cdots,
\end{equation}  where $N=n+1$.

A combination of \eqref{remainder-even-sgn} and
 \eqref{remainder-odd-sgn} gives \eqref{the-remainder-sgn}.  Thus completes the proof of the theorem. \hfil\qed

\section{Lemma \ref{lem-alpha-k-estimate}: Estimating of the coefficients $\alpha_k$}

To prove Lemma \ref{lemma-R(N)}, first we estimate the coefficients $\alpha_k$, or, more precisely,  the quantities  $\tilde \alpha_{2k}= \frac {(-1)^{k+1} \alpha_{2k}} {(2k-1)!}$ and   $\tilde \alpha_{2k+1}= \frac {(-1)^{k+1} \alpha_{2k+1}} {(2k)!}$ for   $k=1,2,\cdots$. These are positive constants; cf. \eqref{alpha-even-sgn} and \eqref{alpha-odd-sgn}.
As a preparation, we give a brief account of the analytic continuation of the hypergeometric function. The reader is referred to \cite[Sec. 3.2]{li-liu-xu-zhao2} for full details.
We denote
\begin{equation}\label{varphi}
 \varphi(z)=F\left(\sin^2 \frac{z}{4}\right )=  F\left(\frac{1}{2},\frac{1}{2};1;\sin^2 \frac{z}{4}\right)~~\mbox{for}~~\Re z \in (-2\pi, 2\pi)\cup (2\pi, 6\pi).
\end{equation}Then the piecewise-defined function
\begin{equation}\label{varphi-continuation}
 v(z)=\left \{\begin{array}{ll}
                \varphi(z),  & 0\leq \Re z < 2\pi, \\
                \varphi(z) \pm 2i \varphi(z-2\pi),     & 2\pi< \Re z < 4\pi~\mbox{and}~\pm\Im z>0
              \end{array}\right .
\end{equation}furnishes an analytic continuation of $\varphi(z)$ in \eqref{varphi} from the strip $\Re z \in [0, 2\pi)$ to the cut strip $0\leq \Re z<4\pi$ and $z\not\in [2\pi, +\infty)$. What is more, we have the connection formula (see \cite[(3.17)]{li-liu-xu-zhao2})
 \begin{equation}\label{v-connection}
v(z)=v_A(z) -\frac 2 \pi \varphi (z-2\pi)\ln \left (2\pi -z  \right )
\end{equation} for $0< \Re z <4\pi$, with  $v_A(z)$ being  analytic in  the strip, and the branch of the  logarithm being  chosen as $\arg (2\pi -z)\in (-\pi, \pi)$.
We proceed to show that

\noindent\begin{lem}\label{lem-alpha-k-estimate} It holds
\begin{equation}\label{alpha-2k-estimate}
\frac{1.9621}{2k-1}\frac{1}{(2\pi)^{2k}}\leq \tilde \alpha_{2k}\leq\frac{2.2032}{2k-1}\frac{1}{(2\pi)^{2k}}
,~~ k=9, 10,11, \cdots
\end{equation}
 and
\begin{equation}\label{alpha-2k+1-estimate}
\frac{4\ln(2k+1)+0.6551}{\pi(2\pi)^{2k+1}}\leq \tilde \alpha_{2k+1}\leq\frac{4\ln(2k+1)+2.2048}{\pi(2\pi)^{2k+1}},~~k=9,10,11,\cdots .
\end{equation}

 \end{lem}

\begin{figure}[t]
\begin{center}
\includegraphics[height=6cm]{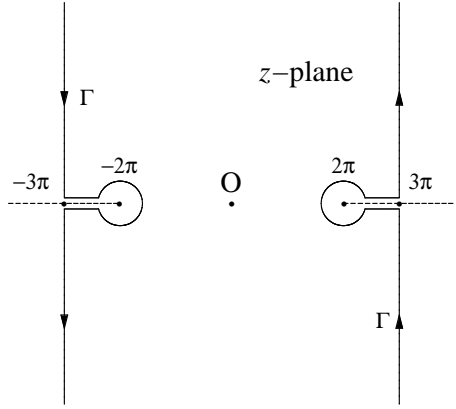}
 \caption{The deformed contour $\Gamma$: the oriented curve (see \cite[Fig.\;2]{li-liu-xu-zhao2}).}
 \label{contour-Gamma}
\end{center}
\end{figure}

\noindent
{\bf {Proof}}.
We understand  \eqref{generating-function-2k1} as a generating relation for $\tilde \alpha_{2k}$. Using the Cauchy integral formula, and in view of \eqref{u-function}, we have
\begin{equation*}
\tilde \alpha_{2k}=\frac{1}{2 \pi i}\oint\frac  {u(z)\cos( {z}/{4})}{z^{2k+1}} dz=\frac{1}{8\pi i}\int_{\Gamma} \frac{v(z)}{\sin({z}/{4})}\frac {dz}{z^{2k}},
\end{equation*}
where initially the  integration path $\Gamma$ is a loop encircling the origin anti-clockwise, and is then deformed to the oriented curve   illustrated in Figure \ref{contour-Gamma}; see also \cite[Fig.\;2]{li-liu-xu-zhao2}, and $v(z)$ is the function defined in \eqref{varphi-continuation}. From \eqref{v-connection}, paying attention to the symmetric properties of  $v(z)$ and $\Gamma$, we have
\begin{equation}\label{alpha-2k-tuta}
\tilde \alpha_{2k}=\frac{1}{4\pi i}\int_{\Gamma_v}\frac{v(z)}{\sin({z}/{4})}\frac { dz }{z^{2k}}
-\frac{1}{4\pi i}\int_{\Gamma_l}\frac{\frac{2}{\pi}\varphi(z-2\pi)\ln(2\pi-z)}{\sin({z}/{4})}\frac {dz}{ z^{2k}}:=I_v+I_l,
\end{equation}where $\Gamma_v$ is the vertical part $\Re z=3\pi$, and $\Gamma_l$ is the remaining right-half part of $\Gamma$, consisting of a circular part around $z=2\pi$, and a pair of horizontal  line segments, respectively along the upper and lower edges of $(2\pi, 3\pi)$,  joining the circle with the vertical line; see Figure \ref{contour-Gamma}.

First, straightforward   calculation gives
\begin{equation}\label{even-Iv}
\left|I_v\right|\leq \frac {M_v}{4\pi} \frac 1 {(3\pi)^{2k}} \int^\infty_{-\infty} \frac {dy}{\big | \sin \frac {3\pi+iy}4\big |}=       \frac{M_v}{2\pi(3\pi)^{2k}}B\left (\frac{1}{4},\frac{1}{4}\right )\approx\frac{3.1153\cdots}{(3\pi)^{2k}},
\end{equation}
where $|v(3\pi+iy)| \leq \sqrt 5  \max_{y\in \mathbb{R}} \left |\varphi(\pi+iy)\right |\leq M_v=2.6393\cdots$; cf. \cite[p.297]{li-liu-xu-zhao2},   and     $B\left (\frac{1}{4},\frac{1}{4}\right )$ is the Beta function.

Now we turn to
  the dominant  part $I_l$. It is readily seen that
  \begin{equation}\label{even-Il}
    I_l=\frac 1 \pi\int^{3\pi}_{2\pi} \frac {\varphi(x-2\pi)}{ \sin (x/4)}\frac {dx}{x^{2k}}= \frac{1}{\pi}\int_{2\pi}^{3\pi}\frac{dx}{x^{2k}} +\frac{1}{\pi}\int_{2\pi}^{3\pi}\left\{\frac{g(x)-1}{x-2\pi}\right\}\frac{(x-2\pi)}{x^{2k}}dx,
  \end{equation}
  where $g(x)=\frac{\varphi(x-2\pi)}{\sin({x}/{4})}$  such that  $g(2\pi)=1$.
 One can see that  $\frac{g(x)-1}{x-2\pi}$ is positive and monotone increasing for $x\in (2\pi, 3\pi]$ since
 \begin{equation*}
 \frac{g(x)-1}{x-2\pi}=\left\{ \frac {\sin(t/4)}{t\cos(t/4)}\right \} \left\{ \frac {\varphi(t) -1}{\sin(t/4)}\right\} +\left\{ \frac 1 t \left (\frac 1 {\cos(t/4)}-1\right)\right\},~~t=x-2\pi,
 \end{equation*} and each right-hand side term in the curly braces is positive and monotone increasing for $t\in (0, \pi]$; see \cite[p.299]{li-liu-xu-zhao2} for the monotonicity of $\frac {\varphi(t) -1}{\sin(t/4)}$. Therefore, we have for $x\in (2\pi, 3\pi]$,
  \begin{equation*}
 0\leq \frac{g(x)-1}{x-2\pi}\leq \frac{g(3\pi)-1}{\pi}:= M_g=\frac 1 \pi \left [\sqrt 2\; F\left (\frac{1}{2},\frac{1}{2};1;\frac{1}{2}\right )-1\right]=  0.2130\cdots.
 \end{equation*}
Substituting it into \eqref{even-Il}, we have
\begin{equation}\label{even-I2}
I_l=\frac{2}{2k-1}\frac{1}{(2\pi)^{2k}}+\frac{\delta_{l,k}}{(2k-1)(2\pi)^{2k}}
\end{equation}
with  $-3\left (\frac{2}{3}\right )^{2k}<\delta_{l,k}\leq\frac{2\pi M_g}{k-1}$.
Further  substituting
\eqref{even-Iv} and \eqref{even-I2} into \eqref{alpha-2k-tuta} gives
\begin{equation}\label{alpha-2k-tuta-approx}
\tilde \alpha_{2k}=\frac{2}{2k-1}\frac{1}{(2\pi)^{2k}}+\frac{\delta_k}{(2k-1)(2\pi)^{2k}}
\end{equation}
with
\begin{equation*}
 -\left \{ 3.1153\;(2k-1)+3\right \} \left (\frac 2 3\right )^{2k} <\delta_{k} <   3.1153\;(2k-1)  \left (\frac 2 3\right )^{2k}  +
\frac{2\pi M_g}{(k-1)}.
\end{equation*}Hence
for $k\geq 9$,  we obtain  the inequalities in  \eqref{alpha-2k-estimate}.\vskip .3cm

Now we turn to the   inequality  \eqref{alpha-2k+1-estimate} for the  odd terms.
From \eqref{u-function} and \eqref{tlide a odd} we have
\begin{equation*}
\tilde \alpha_{2k+1}=\frac{1}{2\pi i}\oint \frac {u(z)\sin({z}/{4}) dz }{z^{2k+2}}=\frac{1}{8\pi i}\int_{\Gamma}\frac{v(z)dz}{\cos({z}/{4}) \; z^{2k+1}},
\end{equation*}where $\Gamma$ is the same path illustrated in Figure \ref{contour-Gamma}.
Then, in view of the connection formula  \eqref{v-connection}, we may write
\begin{equation}\begin{aligned} \label{odd-coefficients-sum}
\tilde \alpha_{2k+1}=&   \frac{1}{4 \pi i}\int_{\Gamma_v}\frac{v(z)}{\cos{\frac{z}{4}}}\frac {dz}  { z^{2k+1}  }+
  \frac{1}{4\pi i}\int_{\Gamma_l}\frac{v_{A}(z)}{\cos{\frac{z}{4}}}\frac {dz}{z^{2k+1}}    \\
&  -\frac{1}{4\pi i}\int_{\Gamma_l}\frac{\frac{2}{\pi}\varphi(z-2\pi)\ln(2\pi-z)}{\cos{\frac{z}{4}}}\frac {dz}{z^{2k+1}}
  :=J_v+J_a+J_l, \end{aligned}\end{equation}
where the integration paths $\Gamma_v$ and $\Gamma_l$ are the same as in \eqref{alpha-2k-tuta}; see Figure \ref{contour-Gamma}.
We note that the   procedure   in  \cite[Sec.3.3]{li-liu-xu-zhao2} applies here, with minor modifications. Case by case estimating gives
\begin{equation}\label{tlide-iv}
\left|J_v\right|\leq \left \{ \frac{M_v}{2  \pi} B\left (\frac{1}{4},\frac{1}{4}\right )\right\}  \frac 1  {(3\pi)^{2k+1}}\approx\frac{3.1153\cdots}{(3\pi)^{2k+1}};
\end{equation}see   \eqref{even-Iv}. Also, picking up the residue at $z=2\pi$ yields
\begin{equation}\label{tlide-ia}
J_a=\frac{2v_A(2\pi)}{(2\pi)^{2k+1}}=\frac{16\ln{2}}{\pi(2\pi)^{2k+1}},
\end{equation}where $v_A(2\pi)=\frac {8 \ln 2} \pi$; see \cite[(3.21)]{li-liu-xu-zhao2}.
The dominant contribution comes from the last integral $J_l$. We follow the steps in \cite[pp.299-301]{li-liu-xu-zhao2}, and eventually obtain
\begin{equation}\label{tlide-il}
J_l=\frac{4\ln(2k+1) -  (4\gamma+4\ln(2\pi))+\delta_{l,k}       }{\pi(2\pi)^{2k+1}},
\end{equation}where
$\left|\delta_{l,k}\right|<2M_\varphi +\frac \pi k \tilde M_f +\frac 4 {k+\frac 1 2}  e^{-k-\frac  1 2} $ for positive integers $k$ with $M_\varphi=\frac {e-1}{2e}$; see \cite[(3.23)]{li-liu-xu-zhao2}, and such  that
\begin{equation*}
  0< \frac {\varphi(x-2\pi)}{\sin\frac {x-2\pi} 4} -\frac 1 {\frac {x-2\pi} 4} =\left \{  \frac {\varphi(t)-1}{\sin\frac {t} 4} \right\} +\left\{ \frac 1  {\sin\frac {t} 4}  -\frac 1 {\frac {t} 4}\right\}
\leq \tilde M_f=\sqrt 2 F\left (\frac{1}{2},\frac{1}{2};1;\frac{1}{2}\right ) -\frac 4\pi  \end{equation*}
for $x\in (2\pi, 3\pi]$, or, $t\in (0, \pi]$ for $t=x-2\pi$. Here use has been made of the fact that both terms in the curly braces are monotone increasing positive functions for $t\in (0, \pi]$; cf. the derivation of \eqref{even-I2}.
Now substituting \eqref{tlide-iv}, \eqref{tlide-ia} and \eqref{tlide-il} into \eqref{odd-coefficients-sum} yields
\begin{equation}\label{alpha-tidle-2k+1}
\tilde \alpha_{2k+1}=\frac{4\ln(2k+1)
 + ( 16\ln 2- 4\gamma-4\ln(2\pi))+\delta_{k}
  }{\pi(2\pi)^{2k+1}},
\end{equation}
where $\left|\delta_{k}\right|\leq   3.1153 \pi \left  (\frac 2 3\right )^{2k+1}+    2M_\varphi +\frac \pi k \tilde M_f +\frac 4 {k+\frac 1 2}  e^{-k-\frac  1 2}$ for positive integers $k$, which is monotone decreasing in $k$.
Straightforward calculation from \eqref{alpha-tidle-2k+1} yields \eqref{alpha-2k+1-estimate} for $k\geq 9$.

Thus, we complete the proof of Lemma \ref{lem-alpha-k-estimate}.     \hfil\qed
\vskip .3cm

For later use, we need the following corollary:

\noindent\begin{cor}\label{cor-ratio}
Assume that $\tilde \alpha_k$ are the positive constants in Lemma \ref{lem-alpha-k-estimate}. Then we have
\begin{equation}\label{ratio-even-even}
\frac {\tilde\alpha_{2k}}{\tilde\alpha_{2k+2}} < \frac {254} 5=50.8~~\mbox{for}~~ k=5,6,\cdots,
\end{equation}
\begin{equation}\label{ratio-odd-odd}
\frac {\tilde\alpha_{2k+1}}{\tilde\alpha_{2k+3}} < 43~~\mbox{for}~~ k=0,1,2,\cdots,
\end{equation}and
\begin{equation}\label{ratio-even-odd}
  (2k+1) \frac {\tilde\alpha_{2k+2}}{\tilde\alpha_{2k+1}} < 0.12~~\mbox{and}~~(2k+1) \frac {\tilde\alpha_{2k+2}}{\tilde\alpha_{2k+3}} <3.7~~  \mbox{for}~~ k=0,1,2,\cdots.
\end{equation}
\hfil\qed
\end{cor}
The results follow accordingly  from Lemma \ref{lem-alpha-k-estimate} and Table \ref{tabel-little-k}.  To obtain \eqref{ratio-odd-odd} one may have to evaluate the ratio $\frac {\tilde\alpha_{2k+1}}{\tilde\alpha_{2k+3}}$ up to $k=13$, such that $\frac {\tilde\alpha_{2k+1}}{\tilde\alpha_{2k+3}}= 38.578, 38.679, 38.762,  38.829,  38.886$ for $k=9, 10, 11, 12, 13$.

\noindent
\begin{table}[h]
  \centering
 \begin{tabular}{|c|c|c|c|c|c|c|c|c|c|}
   \hline
   $k$                                                &  $1$        &   $2$      &  $3$       &    $4$     &  $5$       &  $6$       &  $7$        &  $8$        & $9$     \\
   \hline
  $\frac {\tilde\alpha_{2k}}{\tilde\alpha_{2k+2}} $   &  $56.305$   & $60.184$   & $57.345$   & $53.150$   & $49.797 $  & $47.533 $  &  $46.044 $  & $45.031 $   & $44.303$\\
  \hline
  $\frac {\tilde\alpha_{2k-1}}{\tilde\alpha_{2k+1}} $ &  $21.333$   & $30.720 $  & $34.632 $  & $36.358 $  & $37.227 $  & $37.730 $  &  $38.055 $  & $38.282 $   &$38.450 $\\
  \hline
  $ \frac { (2k-1)\tilde\alpha_{2k}}{\tilde\alpha_{2k-1}}$&      $0.1041$ & $0.1184$  &  $ 0.1007$   &    $0.0851$   &   $0.0749$  &  $0.0684$   &  $0.0642$   &  $0.0612$   &  $0.0590$
  \\  \hline
  $ \frac { (2k-1)\tilde\alpha_{2k}}{\tilde\alpha_{2k+1}}$&      $2.2222$ & $3.6373 $  &  $3.4884 $   &    $3.0964 $   &   $2.7884 $  &  $2.5822 $   &  $2.4432 $   &  $2.3438 $   &  $2.2681 $
  \\
\hline
   \end{tabular}
  \caption{The first few ratios. Calculation conducted using Maple, based on \eqref{coefficient-alpha}-\eqref{d-k-s}.}\label{tabel-little-k}
\end{table}

\section{Proof of Lemma \ref{lemma-R(N)}}

Now that we have proved Lemma \ref{lem-alpha-k-estimate}, we turn to the proof of Lemma \ref{lemma-R(N)}.\vskip .3cm

\noindent
{\bf{Proof of Lemma \ref{lemma-R(N)}}}. To prove \eqref{R-even-sgn}, the idea is as follows: First we show that
\begin{equation}\label{r-2l-2j}
   (-1)^l r_{2l, 2j+2} >0~~\mbox{for}~~j=l+1,l+2,\cdots,~~l=0,1,2,\cdots,
\end{equation}
and
\begin{equation}\label{r-2l-difference}
   (-1)^{l+1} \left (  r_{2l, 2j+1} +\frac 1 2  r_{2l, 2j+2}\right ) >0~~\mbox{for}~~j=l+1,l+2,\cdots,~~l=0,1,2,\cdots.
\end{equation}
Then \eqref{R-even-sgn} follows immediately
from \eqref{r-2l-2j} and \eqref{r-2l-difference}  since $x=1/N\in (0, 1/2]$ for $N\geq 2$, and
\begin{align*}
   \tilde R_{2l}(N) &=  \sum^\infty_{j=l+1} (-1)^{l+1} \left (r_{2l, 2j+1}+   r_{2l, 2j+2}\; x  \right ) x^{2j+1}\\
   &\geq \sum^\infty_{j=l+1} (-1)^{l+1} \left (r_{2l, 2j+1}+  \frac 1 2 r_{2l, 2j+2}  \right ) x^{2j+1}\\
   & >0
 \end{align*}
 for $N=n+1=2, 3, \cdots$, and $l=0,1,2,\cdots$.\vskip .3cm

The above idea is   simple, yet the verification of \eqref{r-2l-2j} and \eqref{r-2l-difference}   is quite complicated. We begin with   \eqref{r-2l-2j}.    First,   a combination of  \eqref{coefficient-alpha} and \eqref{R-l-coefficient} gives
\begin{equation*}
 (-1)^{l}  r_{2l, 2l+4}=d_{2l+2, 2l+4} \left \{(-1)^{l} \alpha_{2l+2} \right \}+ \left\{ -d_{2l+1, 2l+4}\right\} \left\{ (-1)^{l+1} \alpha_{2l+1}\right\} >0
\end{equation*}for $l=0,1,2,\cdots$. Here use has been made of \eqref{alpha-even-sgn}, \eqref{alpha-odd-sgn}, and the facts that  $d_{2l+2, 2l+4}>0$ and $d_{2l+1, 2l+4}<0$. Hence   \eqref{r-2l-2j} is true for $j=l+1$.
Therefore,  we need only to prove \eqref{r-2l-2j} for $j=l+2, l+3, \cdots$. In view of \eqref{R-l-coefficient}, it  suffices  to show,  by an induction argument, that
\begin{equation}\label{r-2l-2j-even-odd}
  r_{l,j}^E:= (-1)^{l+1} \sum^l_{k=0} d_{2k, 2j+2} \alpha_{2k} >0~~\mbox{and}~~  r_{l,j}^O:=(-1)^{l+1} \sum^{l-1}_{k=0} d_{2k+1, 2j+2} \alpha_{2k+1} \geq 0
\end{equation}for $j=l+2,l+3,\cdots$ and $l=0,1,2,\cdots$, where $\alpha_0=-1$.
\vskip .3cm

\noindent
\underline{{\it{Proving $r_{l,j}^E >0$ for $j=l+2,l+3,\cdots$ and $l=0,1,2,\cdots$}}:}

Straightforward verification shows that the first inequality in \eqref{r-2l-2j-even-odd}  holds for $l=0$: We see from   \eqref{R-l-coefficient} and  \eqref{d-0-s}  that
   \begin{equation*}
 r_{0,j}^E=r_{0, 2j+2}=\left ( \frac 1 {j+1}-\frac 1 {2j+1}\right ) +\frac 1 {8j} >0~~\mbox{for}~~j=2,3,\cdots.
   \end{equation*}
Similarly, from \eqref{coefficient-alpha} and \eqref{R-l-coefficient} we have $r_{1,j}^E= \alpha_2 d_{2,2j+2}+\alpha_0 d_{0, 2j+2}$. Hence
 \begin{equation*}
 r_{1,j}^E=  \frac 5 {192} \left [ 2j+2+\frac {2j-1} 4\right ] -  \left [ \frac 1 {j+1}-\frac 1 {2j+1}+\frac 1 {8j}\right ]> \frac {5(j+1)}{96}-\frac 2 {3(j+1)} >0
 \end{equation*}for $j=3,4,\cdots$.  Thus the first inequality in \eqref{r-2l-2j-even-odd} is true  for $l=1$.

Now assume  \eqref{r-2l-2j-even-odd}   for a non-negative integer $l$, then, replacing $l$ with $l+2$, we have
\begin{equation}\label{r-2l-2j-induction}
   r_{l+2,j}^E =  r_{l,j}^E    + \tilde \alpha_{2l+4} \left [ (2l+3)!  d_{2l+4, 2j+2} - \frac { \tilde \alpha_{2l+2}}{ \tilde \alpha_{2l+4} }   (2l+1)!  d_{2l+2, 2j+2}\right ] >0
\end{equation}for $j=l+4,l+5,l+6,\cdots$.
Indeed, if we write
\begin{equation*}
(2l+1)!  d_{2l+2, 2j+2}= \frac {(2j+2l+2) (2j)! } {(2j-2l)!} + \frac {  (2j-1)! } {4(2j-2l-2)!}:=A_{l}+B_{l},
\end{equation*}
Then, noting that for $l\geq 0$ and $j-l\geq 5$, in view of  \eqref{ratio-even-even} and Table \ref{tabel-little-k},  we have
\begin{align*}
    & A_{l+1}+B_{l+1}-\frac{\tilde \alpha_{2l+2}}{ \tilde \alpha_{2l+4}}\left(A_{l}+B_{l}\right)\\
    &\geq (2j-2l-1)(2j-2l) A_l + (2j-2l-3) (2j-2l-2)B_l- 61  \left( A_{l}+B_{l}\right)\\
     &\geq 90A_l+56 B_l- 61 \left(A_{l}+B_{l}\right)\\
     &>0,
\end{align*}
since $A_l> 4 B_l$ by straightforward verification. Alternatively, applying  \eqref{ratio-even-even},    for $l\geq 4$ and $j-l\geq 4$,  we can modify the above inequalities  to give
  \begin{equation*}
 A_{l+1}+B_{l+1}-\frac{\tilde \alpha_{2l+2}}{ \tilde \alpha_{2l+4}}\left(A_{l}+B_{l}\right)
    \geq 56A_l+30 B_l- \frac {254} 5  \left(A_{l}+B_{l}\right)>0
 \end{equation*}
The remaining cases, namely $j=l+4$ with $l=0,1,2,3$, can  be justified by direct calculation: The values  of
\begin{equation*}
(2l+3)!  d_{2l+4, 2l+10} - \frac { \tilde \alpha_{2l+2}}{ \tilde \alpha_{2l+4} }   (2l+1)!  d_{2l+2, 2l+10}=
62.9,\; 1004.5,\; 0.66\times 10^6,\; 0.33\times 10^9,
\end{equation*}respectively for $l=0,1,2,3$. Summarizing  all above, we see the validity of \eqref{r-2l-2j-induction}. Therefore, the first inequality  in \eqref{r-2l-2j-even-odd} is true for $j=l+2,l+3,\cdots$ and $l=0,1,2,\cdots$.
\vskip .3cm

\noindent
\underline{{\it{Proving  $r_{l,j}^O \geq 0$ for $j=l+2,l+3,\cdots$ and $l=0,1,2,\cdots$}}:}

The analysis of $r_{l,j}^O$ is similar to, and simpler than, that of the even terms   $r_{l,j}^E$.  First, for $l=0$, the sum in \eqref{r-2l-2j-even-odd} is empty and thus we understand that  $r_{0,j}^O=0$ for all $j$.
Also, it is readily seen that  $r_{1,j}^O=d_{1, 2j+2} \alpha_1\equiv \frac 3 {16}$ for $j=3,4,\cdots$; cf. \eqref{d-k-s}. Hence, the equality for $r_{l,j}^O$ in \eqref{r-2l-2j-even-odd} also holds  for $l=1$.

Now assume that $r_{l,j}^O\geq 0$ for a non-negative integer $l$ and $j=l+2, l+3,\cdots$. From \eqref{R-l-coefficient}  we may write
\begin{equation}\label{r-2l-2j-odd-induction}
  r_{l+2,j}^O =  r_{l,j}^O-\tilde \alpha_{2l+3}(2l+2)!\;d_{2l+3, 2j+2}
   +\tilde\alpha_{2l+1}(2l)!\;d_{2l+1, 2j+2}:=  r_{l,j}^O+c_+\Delta_{l,j}
\end{equation}with a positive constant $c_+=  \tilde \alpha_{2l+3}  (2l)!\left  | d_{2l+1, 2j+2}\right |$,     and try to prove that
\begin{equation*}
\Delta_{l,j}:=\frac{(2l+2)!\;d_{2l+3, 2j+2}}  {(2l)!\;d_{2l+1, 2j+2}} -\frac{\tilde\alpha_{2l+1}}{\tilde \alpha_{2l+3}}=     \frac { 6j+2l+2 }{ 6j+2l }  (2j-2l-1) (2j-2l) -\frac{\tilde\alpha_{2l+1}}{\tilde \alpha_{2l+3}} >0
\end{equation*}for $j=l+4, l+5,\cdots$.  Here
 the last inequality comes from   \eqref{ratio-odd-odd}. Therefore, from \eqref{r-2l-2j-odd-induction} and  by induction, we have justified the validity of both inequalities in   \eqref{r-2l-2j-even-odd} for
 all $j\geq l+2$ and  $l\geq 0$.  Accordingly,  we  have proved \eqref{r-2l-2j} for  all $j\geq l+1$ and  $l\geq 0$, noting that the only exceptional case $j=l+1$ has been discussed earlier in this section.
\vskip .3cm

\noindent
\underline{{\it{Proving \eqref {r-2l-difference}}}:}

In what follows we proceed to prove \eqref {r-2l-difference}.
First, taking into account the formulas \eqref{coefficient-alpha} and \eqref{R-l-coefficient}, we see that $r_{2l, 2l+3}+ \frac 1 2 r_{2l, 2l+4}$ can be represented as a  linear combination of $\alpha_{2l+1}$ and $\alpha_{2l+2}$. More precisely, substituting in the coefficients $d_{k,s}$; see \eqref{d-0-s} and \eqref{d-k-s}, we have
\begin{equation*}
(-1)^{l+1}\left [ r_{2l, 2l+3}+ \frac 1 2 r_{2l, 2l+4}\right ]=(2l)! \; \tilde\alpha_{2l+1} \left [  \frac {(2l+1)(12l+5)} 8  -\frac {\tilde\alpha_{2l+2}} {\tilde\alpha_{2l+1}} \frac {(2l+2)^2(2l+1)}  2  \right ],
\end{equation*}which is positive for all $l$  since $\frac {\tilde\alpha_{2l+2}} {\tilde\alpha_{2l+1}} < \frac 1 {4(2l+1)}$ for $l\geq 0$; cf. \eqref{ratio-even-odd}. Thus \eqref{r-2l-difference} is true for $j=l+1$,   allowing  us  to just prove \eqref{r-2l-difference} for $j=l+2,l+3,\cdots$ and $l=0,1,2,\cdots$.

For $l=0$, it is readily verified from  \eqref{R-l-coefficient} and   \eqref{d-0-s}   that
\begin{equation*}
-r_{0, 2j+1}- \frac 1 2 r_{0, 2j+2} =\left (\frac 1 {2j}-\frac 1 {2j+2} \right ) + \left ( \frac 1 {2(2j+1)}-\frac 1 {4(2j-1)}-\frac 1 {16j}\right )>0
\end{equation*}for $j=2,3,\cdots$. Here the right-hand side is the sum of positive numbers when $j\geq 3$, and equals to  $\frac {11}{160}$ when $j=2$.
 Hence   \eqref{r-2l-difference} holds for $l=0$.

For $l=1$, recalling that $r_{2,s}=-\frac 5 {192} d_{2,s}+\frac 1 4 d_{1,s} +d_{0,s}$; cf. \eqref{R-l-coefficient}, from \eqref{d-0-s} and \eqref{d-k-s} we may write
 \begin{equation*}
 r_{2, 2j+1}+ \frac 1 2 r_{2, 2j+2} = \frac{ 5 j} {768} +\frac {281}{1536} -  \left (\frac 1 {2j}-\frac 1 {2j+2} \right ) - \left ( \frac 1 {2(2j+1)}-\frac 1 {4(2j-1)}\right )+\frac 1 {16j}.
 \end{equation*}Using   the facts that  $\frac 1 {2j} -\frac 1 { 2j+2 }\leq \frac 1 {8j}$ for   $j\geq 3$,  and  $\frac 1 {2(2j+1)}-\frac 1 {4(2j-1)}<  \frac 1 {8j}$ for   $j\geq 1$,  we have
$r_{2, 2j+1}+ \frac 1 2 r_{2, 2j+2}> \frac{ 5 j} {768} +\frac {281}{1536} -\frac 3 {16j}>0$ for all  $j=3,4,5,\cdots$.  Hence   \eqref{r-2l-difference} holds for $l=1$.

Now assume \eqref {r-2l-difference} for a non-negative integer $l$,  then, from \eqref{R-l-coefficient} we have
\begin{equation*}
(-1)^{l+3}\left (r_{2l+4,2j+1}+\frac{1}{2}r_{2l+4,2j+2}\right )=(-1)^{l+1}\left (r_{2l,2j+1}+\frac{1}{2}r_{2l,2j+2}\right )+O_l+E_l.
\end{equation*}It suffices  to show that $O_l+E_l >0$ for $j=l+4, l+5, \cdots$, where for   $l=0,1,2,\cdots$,
\begin{equation}\label{O-l-def}
 O_l:=(-1)^l \left[\alpha_{2l+3}\left(d_{2l+3, 2j+1}+\frac 1 2 d_{2l+3, 2j+2}\right)+\alpha_{2l+1}\left(d_{2l+1, 2j+1}+\frac 1 2 d_{2l+1, 2j+2}\right)\right]
\end{equation}
  and
\begin{equation}\label{E-l-def}
E_l:=(-1)^l \left[\alpha_{2l+4}\left(d_{2l+4, 2j+1}+\frac 1 2 d_{2l+4, 2j+2}\right)+\alpha_{2l+2}\left(d_{2l+2, 2j+1}+\frac 1 2 d_{2l+2, 2j+2}\right)\right]  .
\end{equation}

We may write
\begin{equation}\label{O-l-representation}
  O_l =\tilde\alpha_{2l+3}(\tilde A_{l+1}+\tilde B_{l+1})-\tilde\alpha_{2l+1}(\tilde A_{l}+\tilde B_{l})= \tilde\alpha_{2l+3}\left [  \tilde A_{l+1}+\tilde B_{l+1} - \frac {\tilde\alpha_{2l+1}} {\tilde\alpha_{2l+3}}
  (\tilde A_{l}+\tilde B_{l}) \right ]
\end{equation} with $\tilde A_l=\frac{(2j-1)!(5j+7l)}{4(2j-2l)!}$ and $\tilde B_l=\frac{(2j-2)!}{4(2j-2l-2)!}$.
Observing that $\tilde A_{l+1} >  (2j-2l)(2j-2l-1) \tilde A_{l}\geq 90 \tilde A_{l}$ and $\tilde B_{l+1} =  (2j-2l-2)(2j-2l-3) \tilde B_{l}\geq 56 \tilde B_{l}$  for $j\geq l+5$ and $l\geq 0$,
and recalling that  $\frac {\tilde\alpha_{2l+1}} {\tilde\alpha_{2l+3}}< 43$ for $l\geq 0$, we have
\begin{equation*}
  O_l \geq   \tilde\alpha_{2l+3}\left [  \tilde A_{l+1}\left (1- \frac {43}{90}\right )    +\tilde B_{l+1} \left (1- \frac {43}{56} \right ) \right ]   \geq   \frac {47}{90}  \tilde\alpha_{2l+3} \tilde A_{l+1},~~j\geq l+5, ~l\geq 0.
\end{equation*}

Now we turn to $E_l$. Similar to the discussion of $O_l$, we may also write
\begin{equation*}
   E_l =\tilde \alpha_{2l+2}   (\tilde C_l-\tilde  D_l  )- \tilde \alpha_{2l+4}  (\tilde C_{l+1}- \tilde D_{l+1} ),
  \end{equation*}where $\tilde C_l=\frac {(2j+1)!}{(2j-2l)!}  -\frac {(2j-1)!}{(2j-2l-2)!}=\frac {(4l+2) (2j-l) (2j-1)!}{(2j-2l)!}$,  and $\tilde D_l= \frac 1 2 \frac {(2j)!} {(2j-2l-1)!} -  \frac 1 8 \frac {(2j-1)!} {(2j-2l-2)!}   -  \frac 1 4 \frac {(2j-2)!} {(2j-2l-3)!}$. It is readily verified that both constants are positive, and such that $  \frac 18  D_l <  \tilde D_l<  \frac 1 2 D_l$ for $j\geq l+4$ and $l\geq 0$  with $D_l=\frac {(2j)!} {(2j-2l-1)!}$. Therefore, we have for $j\geq l+5$ and $l\geq 0$ that
    \begin{equation}\label{O+E-l-inequality}
O_l+E_l >  \frac {47}{90}  \tilde\alpha_{2l+3} \tilde A_{l+1} -  \frac 1 2 \tilde \alpha_{2l+2}  D_l  - \tilde \alpha_{2l+4}  \tilde C_{l+1}: =\frac {(2j-1)!\tilde\alpha_{2l+3}}{(2j-2l-1)!} \Omega ,
 \end{equation}where
  \begin{align*}
 \Omega&=(2j-2l-1)\left [ \frac {47}{360} (5j+7l+7) -\frac  {2 (2l+3)  \tilde \alpha_{2l+4}} {\tilde\alpha_{2l+3}}  (2j-l-1)\right ]-
 \frac  {(2l+1)  \tilde \alpha_{2l+2}} {\tilde\alpha_{2l+3}}\frac j{2l+1}\\
 &\geq 9 \left [ \frac {47}{360} (5j+7l+7) -0.24 (2j-l-1)\right ]-
 \frac {3.7} 3 j \\
 &=\frac {193}{600}j +\frac {2077}{200} l+\frac {2077}{200},
 \end{align*}and thus is positive for $j\geq l+5$ and $l\geq 1$. Here use has been made of \eqref{ratio-even-odd}. For the special case $l=0$ and $j\geq 5$, taking Table \ref{tabel-little-k} into account, again we  have the positivity of $\Omega$:
\begin{equation*}
\Omega \geq (2j-1)\left [\frac {47}{360}(5j+ 7) -  0.24 (2j-1)\right ]-2.23 j=\frac {701}{100}+\frac {6049}{1800}(j-5) +\frac {311}{900}(j-5)^2
\end{equation*}

What remains is the case when $j=l+4$ with $l=0,1,2,\cdots$.  Still we have \eqref{O-l-representation}. Since $\tilde A_{l+1}=\frac {(2j-1)! (12l+27)} { 4 \cdot 6!} > 56 \tilde A_{l}$ and $\tilde B_{l+1}=\frac {(2j-2)!  } { 4 \cdot 4!} = 30 \tilde B_{l}$, in view of  \eqref{ratio-odd-odd} we have
\begin{equation*}
  O_l \geq   \tilde\alpha_{2l+3}\left [  \tilde A_{l+1}\left (1- \frac {43}{56}\right )   -\tilde B_{l+1} \left ( \frac {43}{30} -1 \right ) \right ]   \geq  \left [ \frac {13}{56}-\frac {13}{(2l+7)(12l+27)}\right ]  \tilde\alpha_{2l+3} \tilde A_{l+1},
\end{equation*} from which we see that $O_l>\frac 1 5  \tilde\alpha_{2l+3} \tilde A_{l+1}$ for $j= l+4$ with $l\geq 2$.  As a results, we have a modified version of \eqref{O+E-l-inequality} as $j=l+4$,
\begin{equation*}
O_l+E_l >  \frac {1}{5}  \tilde\alpha_{2l+3} \tilde A_{l+1} -  \frac 1 2 \tilde \alpha_{2l+2}  D_l  - \tilde \alpha_{2l+4}  \tilde C_{l+1}: =\frac {(2j-1)!\tilde\alpha_{2l+3}}{6!} \Omega_4 ,
\end{equation*}where
\begin{equation*}
\Omega_4 = \frac {12l+27} {20}-    \frac  {(2l+3)  \tilde \alpha_{2l+4}} {\tilde\alpha_{2l+3}}  (2l+14) - \frac  {(2l+1)  \tilde \alpha_{2l+2}} {\tilde\alpha_{2l+3}}\frac {l+4} {7(2l+1)}.
\end{equation*}From   \eqref{ratio-even-odd} we readily see that
\begin{equation*}
\Omega_4\geq \frac {12l+27} {20}
-0.12   (2l+14) - {3.7}\times  \frac 1 {7}  = \frac {31}{140} + \frac 9 {25}(l-3),
\end{equation*}  and is positive   for $l\geq 3$.

We fill the last gap  by calculating from \eqref{O-l-def}-\eqref{E-l-def} that
$ O_l+E_l=3.3236, ~ 1.9908, ~  4.3827$, respectively  for  $l=0,1,2$, with $j=l+4$.
 Thus we complete the proof  of \eqref{r-2l-difference},
 and hence of       Lemma \ref{lemma-R(N)}. \hfil\qed

\section{Discussion}

We have proved the conjecture of Granath \cite{Granath}, as stated in Theorem \ref{thm epsilon-N} and Corollary \ref{thm-inequalities}, of which the results of
 Zhao \cite{zhao},  Mortici \cite{Mortici} and Granath \cite{Granath} are special  cases. The asymptotic expansion involved, namely \eqref{Granath-expansion}, corresponds to the special case of Nemes' expansion \eqref{nemes-expansion} in descending powers of $n+h$, with $h=1$.

Earlier in \cite{li-liu-xu-zhao2},   Li {\it et al.} consider the case $h=3/4$; cf. \eqref{beta-alternative} and \eqref{beta-bounds}. According to a result  in  \cite{li-liu-xu-zhao2},  the error due to truncation is bounded in absolute value by, and of the same sign as, the first neglected term for all nonnegative $n$. As an application, we obtain optimal upper and lower bounds up to all orders, holding for all integers $n\geq 0$.

Then, a natural question may arise:
\noindent \begin{qe}(Li, Liu, Xu and Zhao)
 Considering the general expansion in  \eqref{nemes-expansion},    for what $h$ do we have the ``best'' approximation in the sense of \cite[Theorem\;1]{li-liu-xu-zhao2} (or, \eqref{beta-bounds} in the present paper), or in the sense of Theorem  \ref{thm epsilon-N} and Corollary \ref{thm-inequalities}?
 \end{qe}

It is worth noting that the coefficients of the expansion \eqref{nemes-expansion} possess a symmetric property, namely, $g_k(h)=(-1)^k g_k(3/2-h)$. Hence, if we take $h=1/2$,  write $N=n+ 1/2$,  and specify   \eqref{nemes-expansion} as
\begin{equation}\label{expansion-n+1/2}
\pi G_n\sim \ln(16N)+\gamma +\sum_{k=1} ^\infty   \frac{\gamma_k}{N^k}.  \end{equation}
Then it is readily seen that $\gamma_k=(-1)^k\alpha_k$, and hence    $(-1)^{\frac {(l+1)(l+2)} 2} \gamma_{l+1} <0$ for nonnegative integers $l$, very similar to the result in Theorem \ref{coefficient}. Naturally,   analysis similar to what we have conducted in the present paper might lead to
       $
     (-1)^{\frac {(l+1)(l+2)} 2} \epsilon_l(N) <0
 $ for all $n=0,1,2,\cdots$ and $l=0,1,2,\cdots$, where $\epsilon_l(N)$ is the error of \eqref{expansion-n+1/2} due to truncation at the $l$-th order term, with $N=n+ 1/2$.


\begin{thebibliography}{99}
\bibitem{as} M. Abramowitz and I.A. Stegun,
{\it Handbook of Mathematical Functions}, Dover, New York, 1972.


\bibitem{Alzer} H. Alzer, Inequalities for the constants of Landau and Lebesgue, {\it J. Comput. Appl. Math.},      {\bf{139}} (2002),  215-230.

\bibitem{Brutmam} L. Brutman,
A sharp estimate of the Landau constants, {\it J. Approx. Theory.},      {\bf{34}} (1982),  217-220.

\bibitem{chen2014}
C.-P. Chen,  New bounds and asymptotic expansions for the constants of Landau and Lebesgue,  {\it Appl. Math. Comput.},  {\bf 242}  (2014), 790-799.

\bibitem{chen-choi2014}
 C.-P. Chen and  J.  Choi,  Inequalities and asymptotic expansions for the constants of Landau and Lebesgue, {\it Appl. Math. Comput.},   {\bf 248}  (2014), 610-624.


\bibitem{Cvijovic-Klinowski} D. Cvijovi\'{c} and J. Klinowski,
 Inequalities for the Landau constants, {\it Math. Slovaca},      {\bf{50}} (2000), 159-164.


\bibitem{Cvijovic-Srivastava} D. Cvijovi\'{c}  and H.M. Srivastava,
 Asymptotics of the Landau constants and their relationship with hypergeometric functions, {\it Taiwanese J. Math.},      {\bf{13}} (2009), 855-870.



\bibitem{Falaleev} L.P. Falaleev,
 Inequalities for the Landau constants, {\it Sib. Math. J.},      {\bf{32}} (1991), 896-897.

\bibitem{Granath} H. Granath,
On inequalities and asymptotic expansions for the Landau constants, {\it J. Math. Anal. Appl.},      {\bf{386}} (2012), 738-743.


\bibitem{ismail-li-rahman2015}M.E.H. Ismail,  X. Li and M. Rahman, Landau constants and their $q$-analogues,
{\it Anal. Appl.},  {\bf{13}} (2015), 217-231. 



\bibitem{Landau} E. Landau,
Absch\"{a}tzung der koeffizientensumme einer potenzreihe, {\it Arch. Math. Phys.},      {\bf{21}} (1913), 42-50, 250-255.

\bibitem{li-liu-xu-zhao1} Y.-T. Li, S.-Y. Liu, S.-X. Xu and Y.-Q. Zhao, Full asymptotic expansions of the Landau constants via a
difference equation approach, {\it Appl. Math. Comput.},      {\bf{219}} (2012),  988-995.

\bibitem{li-liu-xu-zhao2} Y.-T. Li, S.-Y. Liu, S.-X. Xu and Y.-Q. Zhao, Asymptotics of Landau constants with optimal error bounds, {\it Constr. Approx.},      {\bf{40}} (2014),  281-305.

\bibitem{Mortici} C. Mortici,
Sharp bounds of the Landau constants, {\it Math. Comp},      {\bf{80}} (2011), 1011-1018.

\bibitem{nemes}G. Nemes,
Proofs of two conjectures  on the Landau constants,
 {\it J. Math. Anal. Appl.},  {\bf{388}}  (2012),    838-844.

\bibitem{nemes-nemes}G. Nemes and A. Nemes,
A note on the Landau constants,
 {\it Appl. Math. Comput.},  {\bf{217}}  (2011),    8543-8546.

\bibitem{nist}F. Olver, D. Lozier, R. Boisvert and C. Clark,
{\it NIST handbook of mathematical functions},  Cambridge University Press, Cambridge, 2010.

\bibitem{Watson}G.N. Watson, The constants of Landau and Lebesgue,
 {\it Q. J.  Math. Oxford Ser.},  {\bf{1}}  (1930),    310-318.


\bibitem{Wong-Li1992a} R. Wong and H. Li,
Asymptotic expansions for second-order linear difference equations,
\textit{J. Comput. Appl. Math.}, \textbf{41}
(1992), 65-94.



 \bibitem{zhao}D. Zhao,
Some sharp estimates of the constants of  Landau and Lebesgue,
 {\it J. Math. Anal. Appl.},   {\bf{349}}  (2009),    68-73.


\end{thebibliography}
\end{document}